\documentclass[12pt]{amsart}
\usepackage[centertags]{amsmath}
\usepackage{amsfonts}
\usepackage{amssymb}
\usepackage{graphicx}

\title{Kernels, Truth and Satisfaction}

\author{James H. Schmerl}
\date{\today}

\DeclareMathOperator{\Th}{Th}
\DeclareMathOperator{\Cl}{Cl}
\DeclareMathOperator{\Ht}{ht}
\def\Sent{{\sf Sent}}

\newcommand{\NN}{{\Bbb N}}

\def\pa{{\sf PA}}

\def\NN{{\mathcal N}}
\def\MM{{\mathcal M}}

\begin{document}
\maketitle

\smallskip

 {\bf Prologue.} The well known  Kotlarski-Krajewski-Lachlan Theorem \cite{kkl} says that every model $\MM$  of Peano Arithmetic (\pa) has an elementary extension $\NN \succ \MM$ having a full satisfaction class (or, equivalently, every resplendent model has a full satisfaction class).  Later,   Enayat \& Visser \cite{ev} gave another proof.  According to \cite{ev}, the proof in \cite{kkl}  used some ``rather exotic proof-theoretic technology'', while the proof in \cite{ev}   uses ``a perspicuous method for the construction of full satisfaction classes''. Although not made explicit there, the proof in \cite{ev}, when  stripped to its essentials,  is seen to ultimately depend on showing that certain digraphs have kernels. This is made explicit here. 
 
 There is a lengthy discussion in \S4 of \cite{ev} about the relationship of full satisfaction classes to full truth classes.  Satisfaction classes, which are sets of ordered pairs consisting of a formula in the language of arithmetic and  an assignment for that formula, are exclusively used in \cite{ev}. Truth classes are sets of arithmetic sentences that may also have domain constants. By  \cite[Prop.~4.3]{ev} (whose ``routine but laborious proof is left to the reader"), there is a canonical correspondence between  full truth classes and \mbox{{\em extensional}} full satisfaction classes. The culmination of \cite[\S4]{ev} is the construction of extensional full satisfaction classes.  In \S2 of this paper, we will avoid the intricacies of \cite[\S4]{ev} by working exclusively with  truth classes  to easily obtain the same conclusion.

\bigskip

{\bf \S1.~Digraphs and kernels.} A binary relational structure ${\mathcal A} = (A,E)$ is referred to here as a {\bf directed graph}, or {\bf digraph} for short.{\footnote{Henceforth, ${\mathcal A}$ always denotes a digraph  $(A,E)$. If $B \subseteq A$, then we often identify $B$ with the the induced subdigraph  ${\mathcal B} = (B, E \cap B^2)$.}}  A subset $K \subseteq A$ is a {\bf kernel} 
of ${\mathcal A}$ if for every $a \in A$, $a \in K$ iff whenever  $aEb$, then $b \not\in K$. According to \cite{book}, kernels were introduced by von Neumann \cite{von} and have subsequently found many applications.  For $n < \omega$,  define the binary relation $E^n$ on $A$ by recursion: $xE^0y$ iff $x = y$; $xE^{n+1}y$ iff $xEz$ and $zE^n y$ for some $z \in A$.  A digraph ${\mathcal A}$ is a {\bf directed acyclic graph} ({\bf DAG}) if  whenever  $n < \omega$ and $aE^na$, then $n = 0$.   Some DAGs have kernels while others do not. For example, if $<$ is a linear order of $A$ with no maximum element, then $(A,<)$ is a DAG with no kernel. However, every {\it finite} DAG has a (unique) kernel, as was  first noted in \cite{von}.  

An element $b \in A$ for which there is no $c \in A$ such that $bEc$ is a {\bf sink} of ${\mathcal A}$. We say that ${\mathcal A}$ is {\bf well-founded} if  every nonempty subdigraph of ${\mathcal A}$ has a sink.  Every finite DAG is well-founded, and every well-founded digraph is a DAG having a kernel.  The next proposition, for which we need some more  definitions, says even more is true.
 A subset $D$ of a digraph ${\mathcal A}$  is {\bf closed} if whenever $d \in D$ and $dEa$, then $a \in D$. If $X \subseteq A$ and $k < \omega$, then define $\Cl^{\mathcal A}_k(X)$ by recursion: 
 $\Cl^{\mathcal A}_0(X)  = X$ and 
 $\Cl^{\mathcal A}_{k+1}(X) = X \cup \{a \in A : d Ea$ for some $d \in \Cl^{\mathcal A}_k(X)\}$.  
 Let $\Cl^{\mathcal A}(X) = \bigcup_{k < \omega}\Cl^{\mathcal A}_ k(X)$, which is the smallest closed superset of $X$.

\bigskip

{\sc Proposition 1}: {\em Suppose that ${\mathcal A}$ is a digraph, $D \subseteq A$ is closed, $K_0 \subseteq  D$ is a kernel of $D$, and $A \backslash D$ is well-founded. Then ${\mathcal A}$ has a $($unique$)$ kernel $K$ such that $K_0 = K \cap D$.} 

\bigskip

{\it Proof.} By recursion on ordinals $\alpha$, define $D_\alpha$ so that $D_0 = D$, $D_{\alpha+1} = D_\alpha \cup \{b \in D : b$ is a sink of $A \backslash D_\alpha\}$, and $D_\alpha = \bigcup_{\beta < \alpha} D_\beta$ if $\alpha$ is a limit ordinal. Then, there is $\gamma$ such that $A = D_\gamma$.  For each $\alpha$, there is a unique kernel $K_\alpha$ of $D_\alpha$ such that $K_\beta = K_\alpha \cap D_\beta$ whenever $\beta < \alpha$. Let $K = K_\gamma$. \qed

\bigskip

Let ${\mathcal A}$ be a digraph. If there is $k < \omega$ for which there are 
no $a, b \in A$ such that $aE^{k+1}b$, then ${\mathcal A}$ has {\bf finite height}, and we let $\Ht({\mathcal A})$, the {\bf height} of ${\mathcal A}$, be the least such $k$. If ${\mathcal A}$ has finite height, then it is well-founded.
We say that ${\mathcal A}$ has {\bf local finite height} if for every $m < \omega$ there is $k < \omega$ 
  such that $\Ht(\Cl^{\mathcal A}_m(F)) \leq k$ for every $F \subseteq   A$ having cardinality at most $m$.  
  If ${\mathcal A}$  has local finite height, then it is a DAG.
 Having local finite height is a first-order property: 
 if ${\mathcal B} \equiv {\mathcal A}$ and ${\mathcal A}$ has local finite height, then so does ${\mathcal B}$.

 \bigskip

{\sc Theorem 2}: {\em Every digraph ${\mathcal A}$  having local finite height has an elementary extension ${\mathcal B} \succ {\mathcal A}$ that has a kernel.}

\bigskip

{\it Proof.} This proof is modeled after  Theorem~3.2(b)'s in  \cite{ev}.

To get ${\mathcal B}$ with a kernel $K$, we let $B_0 = \varnothing$, and then obtain  an elementary chain 
${\mathcal A}= {\mathcal B}_1 \prec {\mathcal B}_2 \prec {\mathcal B}_3 \prec \cdots$ and an 
increasing sequence $\varnothing  = K_0 \subseteq K_1 \subseteq K_2 \subseteq \cdots$ such that for every $n < \omega$, $K_n$ is a kernel of  $\Cl^{{\mathcal B}_{n+1}}(B_n)$ . Having these sequences, we let ${\mathcal B} = \bigcup_{n<\omega}{\mathcal B}_{n+1}$ and $K = \bigcup_{n<\omega}K_n$, so that ${\mathcal B} \succ {\mathcal A}$ and $K$   is a kernel of~${\mathcal B}$.  
The  next lemma allows us to get ${\mathcal B}_{n+2}$ and $K_{n+1}$ when we already have 
$B_n$, ${\mathcal B}_{n+1}$ and~$K_n$.

\bigskip

{\sc Lemma 3}: {\em Suppose that ${\mathcal B}_{n+1}$  is a  digraph  having local finite height, 
 $D$ is a closed subset of $B_{n+1}$, and $K_n$ is a kernel of $D$.
Then there are  ${\mathcal B}_{n+2} \succ {\mathcal B}_{n+1}$ and a kernel $K_{n+1}$ of $\Cl^{{\mathcal B}_{n+2}}(B_{n+1})$ such that  $K_n =  K_{n+1} \cap D$.}

\bigskip

To prove Lemma~3, let $\Sigma$ be the union of the following three sets of sentences:

\begin{itemize}

\item $\Th(({\mathcal B}_{n+1}, a)_{a \in B_{n+1}})$;

\item $\{\sigma_{F,k}: k < \omega$ and $F \subseteq B_{n+1}$ is finite$\}$, where $\sigma_{F,k}$ is the sentence  
$$\qquad \forall x \in \Cl_k(F)[U(x) \leftrightarrow \forall y \in \Cl_{k+1}(F) \big(xEy \rightarrow \neg U(y)\big)];$$

\item $\{U(d) : d \in K_n\} \cup \{\neg U(d) : d \in D \backslash K_n\}$.

\end{itemize}
This $\Sigma$ is a set of ${\mathcal L}$-sentences, where ${\mathcal L}  = \{E,U\} \cup B_{n+1}$  and $U$ is a new unary relation symbol.

It  suffices to show that $\Sigma$ is consistent, for then we can let $({\mathcal B}_{n+2},U) \models \Sigma$ and let $K_{n+1} = U \cap \Cl^{{\mathcal B}_{n+2}}(B_{n+1})$.  To do so, we need only show that every  finite subset of $\Sigma$ is consistent. 

Let $\Sigma_0 \subseteq \Sigma$ be finite. Let $k_0 < \omega$ and finite $F_0 \subseteq B_{n+1}$ be  such that if $\sigma_{F,k} \in \Sigma_0$, then $k < k_0$ and $F \subseteq F_0$. Let $D = \Cl_{k_0}^{{\mathcal B}_{n+1}}(F_0)$. Since ${\mathcal B}_{n+1}$ has local finite height, then $D$ has finite height and, therefore, is well-founded.   
 By Proposition~1, ${\mathcal B}_{n+1}$ has a kernel $U$ such that $K_n = U \cap D$. Then, $({\mathcal B}_{n+1},U) \models \Sigma_0$, so $\Sigma_0$ is consistent, thereby proving Lemma~3 
 and also Theorem~2. \qed

\bigskip 

{\sc Corollary 4}: {\em Every  resplendent $($or countable, recursively saturated \!$)$ digraph that  
has local finite height has a kernel.}

\bigskip

{\it Proof}. This is just a definitional consequence of Theorem~2 and the fact \cite{bs} that countable, recursively structures are resplendent. \qed

\bigskip

{\bf \S2.~Truth Classes.}  There are various ways that syntax for arithmetic can defined in a model $\MM$ of $\pa$. It usually makes little difference how it is done, so we will choose a way that is very convenient.

 We will  formalize the language of arithmetic by using just two ternary relation symbols: one for addition and one for multiplication. Suppose that $\MM \models \pa$. For each $a \in M$, we have a constant symbol $c_a$.  Then let ${\mathcal L}^M$ consist of the two ternary relations and all the $c_a$'s. The only propositional connective we will use is the NOR connective $\downarrow$, where  $\sigma_0 \! \downarrow \sigma_1$ is $\neg (\sigma_0 \vee \sigma_1)$. The only  quantifier we will use is the  ``there are none such that'' quantifier $\reflectbox{${\mathsf N}$}$, 
where \reflectbox{${\mathsf N}$}$v\varphi(v)$ is $\forall v[\neg\varphi(v)]$. 
Let {\sf Sent}$^\MM$ be the set of ${\mathcal L}^\MM$-sentences as defined in $\MM$. 
A subset $S \subseteq \Sent^\MM$ is a {\bf full truth class} for $\MM$ provided the following 
hold for every $\sigma \in \Sent^\MM$:

\begin{itemize} 

\item if $\sigma = \sigma_0 \downarrow \sigma_1$, then $\sigma \in S$ iff $\sigma_0, \sigma_1 \not\in S$;

\item If $\sigma = $ \reflectbox{${\mathsf N}$}$v \varphi(v)$, then $\sigma \in S$ iff there is no $a \in M$ such that $\varphi(c_a) \in S$;

\item If $\sigma$ is atomic, then $\sigma \in S$ iff $\MM \models \sigma$.

\end{itemize}

Let $A^\MM = \{\sigma \in \Sent^\MM : $ if $\sigma$ is atomic, then $\MM \models \sigma\}$. Define 
the binary relation $E^\MM$ on $A^\MM$ so that if 
 $\sigma_1, \sigma_2 \in A^\MM$, then $\sigma_2E^\MM \sigma_1$ iff one of the following holds:

\begin{itemize}

\item there is $\sigma_0$ such that $\sigma_2 = \sigma_0 \downarrow \sigma_1$ or $\sigma_2 = \sigma_1\downarrow \sigma_0$;

\item  $\sigma_2 =$ \reflectbox{${\mathsf N}$}$v\varphi(v)$ and $\sigma_1 = \varphi(c_a)$ for some $a \in M$.

\end{itemize}

Let ${\mathcal A} = {\mathcal A}^\MM = (A^\MM, E^\MM)$. Obviously, ${\mathcal A}$ is  a DAG. Moreover,  it  has local finite height: if $F \subseteq A^\MM$ is finite and $m < \omega$, then 
$\Ht(\Cl_m^{\mathcal A}(F)) \leq (2^{m+1} - 1)|F|$. We easily see that $S$ is a full truth class for $\MM$ iff $S$ is a kernel of ${\mathcal A}$. 

We can now infer the following version of the KKL Theorem.

\bigskip

{\sc Corollary 5}: {\em Every resplendent $($or countable, recursively saturated \!$)$ $\MM \models \pa$ has a full truth class.} 

\bigskip

{\it Proof}. Since ${\mathcal A}^\MM$ is definable in $\MM$ and $\MM$ is resplendent,  then ${\mathcal A}^\MM$ is also resplendent. Thus, by Corollary~4, ${\mathcal A}^\MM$ has a kernel, which we have seen is a full truth class for $\MM$. \qed

\bigskip

Corollary~5 can be improved by  replacing $\pa$ with any of its subtheories in which enough syntax is definable.

\bibliographystyle{plain}

\end{document}